\documentclass[10pt]{amsart}
\usepackage{mathptmx}
\usepackage{amsmath}
\usepackage{amssymb}
\usepackage{array}
\usepackage{geometry}
\usepackage[bookmarks=true,colorlinks=true, pdfstartview=FitV, linkcolor=black, citecolor=blue, urlcolor=black]{hyperref}

\usepackage{color}
\definecolor{DarkRed}{rgb}{0.55,.00,0.2}
\definecolor{DarkGrey}{rgb}{0.35,.35,0.35}

\theoremstyle{definition}

\theoremstyle{remark}

\numberwithin{equation}{section}

\begin{document}

\title{Polynomial problems of the  Casas-Alvero type}

\author{S. Yakubovich}
\address{Department of Mathematics, Fac. Sciences of University of Porto,Rua do Campo Alegre,  687; 4169-007 Porto (Portugal)}
\email{ syakubov@fc.up.pt}


\date{\today}


\maketitle

\markboth{\rm \centerline{ S.  YAKUBOVICH}}{}
\markright{\rm \centerline{POLYNOMIAL PROBLEMS }}

\begin{abstract}  We establish necessary and sufficient conditions for an arbitrary polynomial of degree $n$, especially with only real roots,   to be trivial, i.e. to have the form $a(x-b)^n$. To do this,  we derive  new properties of  polynomials and their roots. In particular, it concerns new bounds and genetic sum representations of  the Abel -Goncharov
interpolation polynomials. Moreover, we prove  the  Sz.-Nagy type identities, the  Laguerre and Obreshkov-Chebotarev type inequalities for roots of polynomials  and their derivatives.  As applications these results are associated with the known problem, conjectured by Casas- Alvero in 2001, which says, that  any complex univariate polynomial, having a common root with each of its non-constant derivative  must be a power of a linear polynomial. We investigate
particular cases of the problem,  when the conjecture holds true or, possibly, is  false.
\end{abstract}

\maketitle

\section{Introduction and preliminary results}

It is well known from elementary calculus that an arbitrary
polynomial $f$ with complex coefficients (complex polynomial) of
degree $n \in \mathbb{N}$
$$f(z)= a_0 z^{n}+ a_1 z^{n-1}+ \dots + a_{n-1} z + a_n, \ a_0\neq
0,\eqno(1)$$
having a root $\lambda \in \mathbb{C}$ of multiplicity $\mu, \ 1\le
\mu\le n$,  shares it with each of its derivatives  up to order $\mu-1$,
but $f^{(\mu)}(\lambda)\neq 0$. When $\lambda$ is a unique root of
$f$, it has the  form  $f(z)=a(z-\lambda)^n$, $\mu=n$ and $\lambda$ is
the same root of each derivative of $f$ up to order $n-1$. We will
call such a polynomial  a trivial polynomial.  Obviously, as it
follows from the fundamental theorem of algebra, $f$ has at least two
distinct roots, i.e. a polynomial of degree $n$ is non-trivial, if
and only if its maximum multiplicity of roots $r$ does not exceed
$n-1$.

In 2001 Casas- Alvero \cite{CA} conjectured that an arbitrary polynomial $f$  degree $n
\ge 1$ with complex coefficients  is of the form   $f(z)=
a(z-b)^n, a, b \in \mathbb{C}$,    if   and \ only   if $f$  shares a root with each of its derivatives $f^{(1)}, f^{(2)},
\dots, f^{(n-1)}.$

We will call a possible non-trivial polynomial,  which has a common root with each of its non-constant derivatives   a  CA-polynomial. The conjecture says that there  exist no CA-polynomials.  The problem is still open. However, it is proved for small degrees, for infinitely many degrees, for instance, for all powers $n$,   when $n$ is a prime   (see in \cite{Drai}, \cite{Graf}, \cite{Pols} ).  We observe that such a kind of  CA-polynomial of degree $n \ge 2$  cannot have all distinct roots since at least one root is common   with its  first derivative.  Therefore it has a multiplicity at least 2 and a maximum of possible distinct  roots is $n-1$.

 Our main goal here is to derive necessary and sufficient conditions for an arbitrary polynomial (1) to be trivial.  For example, solving a simple differential equation of the first order, we easily prove that a polynomial is trivial, if and only if it is divisible by its first derivative.  In the sequel we establish other criteria, which will guarantee that an arbitrary polynomial has a unique joint root.

 Without loss of generality one can assume in the sequel that $f$ is a monic polynomial of degree $n$, i.e. $a_0=1$ in (1).  Generally, it has $k$ distinct roots $\lambda_j$  of multiplicities $r_j, \   j= 1,  \dots,  k, 1\le k\le n$ such that
$$ r_1+ r_2+ \dots r_{k} = n\eqno(2).$$
 By $r$ we will denote the maximum of multiplicities (2), $r= \hbox{max}_{1\le j\le k}  (r_j)$,  $r_0= \hbox{min}_{1\le j\le k}  (r_j)$ and by  $  \xi^{(m)}_\nu, \  \nu = 1,\dots, n-m$  the zeros of the $m$-th derivative $f^{(m)}, \ m=1,\dots, n-1.$
 For further needs we specify zeros of the $n-1$st and $n-2$nd derivatives, denoting them by $\xi^{(n-1)}_1=z_{n-1}$ and $\xi^{(n-2)}_2=z_{n-2}$, respectively. It is easy to find another zero of the $n-2$-th derivative, which is equal to $\xi^{(n-2)}_1= 2z_{n-1} -z_{n-2} $.  When  zeros $ z_{n-1},\  z_{n-2}$ are real we write, correspondingly,   $ x_{n-1},\  x_{n-2}$.    The value $z_{n-1}$ is called the centroid. It is a center of gravity of roots and by Gauss-Lucas theorem it is contained in the convex hull of all non-constant polynomial derivatives (see details in \cite{Rah}).

 The paper is structured as follows: In Section 2 we study properties of the Abel-Goncharov interpolation polynomials, including integral and series representations and upper bounds.  Section 3 deals with the Sz.-Nagy type identities and Obreshkov-Chebotarev type inequalities for roots of polynomials and   their derivatives. As applications  new criteria are found for an arbitrary polynomial with only real roots to be trivial.   Section 4 is devoted to the Laguerre type inequalities for polynomials with only real roots to localize their zeros. The final Section 5 contains applications of these results towards  solution of the  Casas-Alvero  conjecture and its particular cases.

 \section{Abel-Goncharov polynomials, their upper bounds and  integral and genetic sum's representations}

We begin, choosing a sequence of complex numbers (repeated terms are permitted)
$z_0, z_1, z_2, \dots, z_{n-1}, n \in \mathbb{N}$, where $z_0 \in \{\lambda_1, \lambda_2, \dots, \lambda_k\},\
z_m \in  \{\xi^{(m)}_1,\  \xi^{(m)}_2, \dots, \xi^{(m)}_{n-m}\}, \  m=1,2,\dots, n-1$, satisfying conditions $f^{(m)}(z_m) =0, \
m =0, 1,\dots,  \ n-1$ and, clearly $f^{(n)} (z)= n!$.     Then we represent  $f(z)$ in the form

$$
    f(z)= z^n + P_{n-1} (z),\eqno(3)
$$
where $P_{n-1} (z)$ is a polynomial of degree at most $n-1$. To determine  $P_{n-1} (z)$  we  differentiate the latter equality $m$ times,  and we  calculate the corresponding derivatives in $z_m$ to obtain
$$
 P_{n-1}^{(m)} (z_m) =  -  \frac{n!}{ (n-m)!}  z_m ^{n- m}  ,   \quad  m= 0,1, \dots, n-1.\eqno(4)
$$
But this is the known Abel-Goncharov interpolation problem (see \cite{Evgrafov}) and the polynomial $P_{n-1}(z)$ can be uniquely determined via the linear system (4) of $n$ equations with $n$ unknowns and triangular
matrix with non-zero determinant.   So, following \cite{Evgrafov}, we derive

$$
     P_{n-1} (z) = -  \sum_{k=0}^{n-1}  \frac{n!}{ (n- k)!}  z_k ^{n- k} G_k(z) , \eqno(5)
$$
where $G_k(z),  k=0, 1,\dots, n-1$ is the system of the Abel-Goncharov polynomials
 \cite{Evgrafov}, \cite{Levinson1},  \cite{Levinson2}.    On the other hand it is known that

$$
     G_n(z)= z^n  -  \sum_{k=0}^{n-1}  \frac{n!}{ (n- k)!}  z_k ^{n- k}  G_k(z).\eqno(6)
$$
Thus comparing with (3),  we  find that
$$
G_n(z)\equiv  G_n\left(z, z_0, z_1, z_2,\dots, z_{n-1}\right)=f(z),
$$
and
$$
 G_n\left(\lambda_j, z_0, z_1, z_2,\dots, z_{n-1}\right) = f(\lambda_j)= 0,   \quad
j= 1,2, \dots,  k.
$$
Plainly, one can relate  possible CA-polynomials with the corresponding Abel-Goncharov polynomials, fixing a sequence $\{z_m\}_0^{n-1}$ such that
$$z_m \in \{\lambda_1, \lambda_2, \dots, \lambda_k\}, \  m=0,1,\dots,\ n-1.$$

Further, It is known \cite{Evgrafov} that the Abel-Goncharov polynomial can
be represented as a multiple integral in the complex plane
$$
G_n(z)=   n! \int_{z_0}^z  \int_{z_1}^{s_{1}} \dots
\int_{z_{n-1}}^{s_{n-1}} d s_{n} \dots d s_{1}.\eqno (6)
$$
Moreover, making  simple changes  of variables in (6),  it can be verified that $G_n(z)$ is
 a homogeneous function of degree $n$ (cf. \cite{Levinson1}). Therefore
$$G_n(\alpha z)= G_n\left(\alpha z, \alpha z_0, \alpha z_1,\dots, \alpha z_{n-1}\right)
 = \alpha^n G_n(z), \ \alpha \neq 0.\eqno(7)$$
The following Goncharov upper bound holds for  $G_n$ (see \cite{Gon}, \cite{Evgrafov}, \cite{Levinson1}, \cite{Ibra})
$$
\left| G_n(z) \right|  \le  \left( |z-z_0| + \sum_{s=0}^{n-2} \left
|z_{s+1}- z_s\right| \right)^n.\eqno(8)
$$
  Let us represent the Abel-Goncharov polynomials $G_n(z)$ in a different  way.  To do this, we will use the following representation of the Gauss hypergeometric function given by relation (2.2.6.1) in \cite{Prud}, namely
$$\int_a^b (z-a)^{\alpha-1} (b-z)^{\beta-1}(z + c)^\gamma dz =  (b-a)^{\alpha+\beta-1} (a+c)^\gamma
B(\alpha,\beta) {}_2F_1 \left(\alpha, -\gamma; \alpha +\beta;
\frac{a-b}{ a+c} \right),\eqno(9)$$ where $\alpha, \beta, \gamma $
are positive integers,  $a, b, c \in \mathbb{C}$ and
$B(\alpha,\beta)$ is the Euler beta-function. So, our goal  will be a
representation of the Abel-Goncharov polynomials  in terms of the
so-called genetic sums considered, for instance, in \cite{Apteka}.
Moreover, this will lead  us to another than (8)  upper bound for these  polynomials.  Indeed,   $G_1(z)= z-z_0$. When $n \ge 2$, we use the   multiple integral  representation (6),  and appealing to  the representation (9), we obtain recursively
$$
G_n(z)=   n! \int_{z_0}^z  \int_{z_1}^{s_{1}} \dots \int_{z_{n-2}}^{s_{n-2}} \  (s_{n-1} - z_{n-1})
 d s_{n-1} \dots d s_{1}
$$
$$=   n! (z_{n-2} - z_{n-1})  \int_{z_0}^z  \int_{z_1}^{s_{1}} \dots \int_{z_{n-3}}^{s_{n-3}} \  (s_{n-2} - z_{n-2})
 {}_2F_1 \left(1, -1 ;  \ 2; \   \frac{ z_{n-2} - s_{n-2}}{ z_{n-2} -  z_{n-1} } \right) d s_{n-2} \dots d s_{1} $$
$$=  n!   \sum_{j_1=0} ^ 1\frac{(-1)_{j_1}(-1)^{j_1} }{(2)_{j_1}}  (z_{n-2} - z_{n-1})^{1-j_1} \int_{z_0}^z
\int_{z_1}^{s_{1}} \dots \int_{z_{n-3}}^{s_{n-3}} \  (s_{n-2} - z_{n-2})^{1+j_1}  d s_{n-2} \dots d s_{1} .$$
Hence, employing  properties of the Pochhammer symbol and repeating this process, we find
$$
G_n(z)=  n!   \sum_{j_1=0} ^ 1\frac{(z_{n-2} - z_{n-1})^{1-j_1}  }{(2)_{j_1} (1-j_1)!}   \int_{z_0}^z  \int_{z_1}^{s_{1}}
\dots \int_{z_{n-3}}^{s_{n-3}} \  (s_{n-2} - z_{n-2})^{1+j_1}  d s_{n-2} \dots d s_{1} $$
$$= n!   \sum_{j_1=0} ^ 1  \sum_{j_2=0} ^{1+j_1} \frac{(z_{n-2} - z_{n-1})^{1-j_1} (z_{n-3} - z_{n-2})^{1+j_1-j_2}  }
{(2)_{j_2} (1-j_1)!(1+j_1 -j_2)!}   \int_{z_0}^z  \int_{z_1}^{s_{1}} \dots \int_{z_{n-4}}^{s_{n-4}} \
(s_{n-3} - z_{n-3})^{1+j_2}  d s_{n-3} \dots d s_{1} .$$
Continuing to calculate  iterated integrals with the use of (9), we
arrive  finally  at the following genetic sum  representation of
the Abel-Goncharov polynomials ($j_0 = j_n=0,\   z_{-1}\equiv z $)
$$
G_n(z)=   n!   \sum_{j_1=0} ^ 1  \sum_{j_2=0} ^{1+j_1} \dots  \sum_{j_{n-1}=0} ^{1+j_{n-2}}  \
 \prod_{s=0}^{n-1} \frac{ (z_{n-2-s}- z_{n-1-s} )^ {1+ j_s- j_{s+1}} }{ (1+ j_s- j_{s+1})!}.\eqno(10) $$
Analogously, we derive the genetic  sum representation for  the $m$-th derivative
$G_{n}^{(m)}(z)$, namely  ($j_0  =0 $)
 $$
G_n^{(m)} (z)=   n!   \sum_{j_1=0} ^ 1  \sum_{j_2=0} ^{1+j_1} \dots  \sum_{j_{n-1-m }=0} ^{1 +j_{n-2-m}}
  \frac{ (z - z_{m} )^ {1 + j_{n-1-m}} }{ (1 + j_{n-1-m})!} \
 \prod_{s=0}^{n-2-m} \frac{ (z_{n-2-s}- z_{n-1-s} )^ {1 + j_s- j_{s+1}} }{ (1+ j_s- j_{s+1} )!},\eqno(11) $$
where $m= 0,1,\dots, n-1$.

Meanwhile,  the Taylor expansions of $ G_n^{(m)} (z)$ in the neighborhood of points  $z_m$ give the formulas
$$G_n^{(m)} (z)  = \frac{n!}{(n- m)!} (z-z_m)^{n-m}  + \frac {G_n^{(n-1)} (z_m)}{(n-m-1)!} (z-z_m)^{n-m- 1} + \dots
+  G_n^{(1+m)} (z_m)  (z-z_m),\eqno(12)$$
where $m= 0,1,\dots, n-1$.  Thus comparing coefficients in front of $(z-z_m)^s, \  s= 1, \dots , n- m-1$ in (11) and (12), we find the values of derivatives $G_n^{(s+ m)} (z_m)$ in terms  of $z_m, z_{m+1}, \dots, z_{n-1}$. Precisely, we obtain   ($j_0  =0 $)
$$G_n^{(s+m)} (z_m)   =   n!   \sum_{j_1=0} ^ 1  \sum_{j_2=0} ^{1+j_1} \dots
\sum_{j_{n-2-m }=0} ^{1 +j_{n-3-m}}   \frac{ (z_m - z_{m+1} )^ {2 + j_{n-2-m}-s} }{ (2 + j_{n- 2-m}-s)!} \
 \prod_{l=0}^{n-3-m} \frac{ (z_{n-2-l}- z_{n-1-l} )^ {1 + j_l- j_{l+1}} }{ (1+ j_l- j_{l+1} )!},\eqno(13)$$
where $s=1,2,\dots, n-m, \  m=0,1,\dots, n-1$.

Finally, in this section,  we will establish an   upper bound for the Abel-Goncharov polynomials (cf. (8)).    We have

{\bf Theorem 1}. {\it Let $z, z_0, z_1, z_2,\dots, z_{n-1} \in
\mathbb{C},\ n \ge 1$. The following  upper bound holds for the
Abel-Goncharov polynomials
$$| G_n\left(z, z_0, z_1, z_2,\dots, z_{n-1}\right)| \le
\sum_{k_0=0}^1  \sum_{k_1=0} ^{2-k_0} \dots  \sum_{k_{n-2}=0} ^{n-1
- k_0-k_1-\dots- k_{n-3}}  \ {n \choose k_0,  k_1,  \dots,  k_{n-2},  \
n - k_0-k_1-\dots- k_{n-2} }$$
$$\times |z-z_0|^{n - k_0-k_1-\dots- k_{n-2}}  \prod_{s=0}^{n-2} |z_{n-2-s}- z_{n-1-s}|^ {k_s},\eqno(14) $$
where $ z_{-1}\equiv z$ and
$${n \choose l_0,  l_1,  \dots,  l_m} = \frac{n!}{l_0! l_1! \dots
l_m!}, \ l_0+l_1+\dots +l_m= n$$ are multinomial coefficients. }

\begin{proof} In fact, making simple substitutions $k_s= 1+j_s- j_{s+1}, \ s=0,1,\dots, n-1,  j_0=j_n=0$
and writing identity (10) for the Abel-Goncharov polynomials (6), we estimate their  absolute value, coming out immediately with inequality (14). 

\end{proof}

\section{Sz.-Nagy type identities for roots of polynomials and their derivatives}

In this section we prove  Sz.-Nagy type identities \cite{Rah} for zeros of monic polynomials with complex coefficients and their derivatives.  All   notations of  roots and their multiplicities given in Section 1 are involved.

    We begin with

{\bf  Lemma 1.}  {\it Let $f$ be  a monic polynomial of degree $n \ge 2$ with complex coefficients,  $m= 0,1,\dots, n-2$ and  $z \in \mathbb{C}$. Then the following Sz.-Nagy type identities, which  relate the roots of $f$ and its $m$-th derivative,  hold }

$$z_{n-1} - z  = {1\over n} \sum_{j=1}^k r_{j}(\lambda_j- z) =   {1\over n-m} \sum_{j=1}^{n-m}  (\xi^{(m)}_j - z),\eqno(15)$$

$$   (z_{n-1} - z_{n-2})^2={1\over n(n-1)} \left[ \sum_{j=1}^k r_{j}(\lambda_j- z)^2-  n (z_{n-1} -  z)^2\right]
=  {1\over (n-m)(n-m -1)}$$$$\times  \left[ \sum_{j=1}^{n-m}
(\xi^{(m)}_j - z)^2 -  (n-m) (z_{n-1} -  z)^2\right],\eqno(16)$$

$$   (z_{n-1} - z_{n-2})^2={1\over n^2(n-1)} \sum_{1\le j <  s\le k} r_{j}r_{s}(\lambda_j- \lambda_s)^2
=  {1\over (n-m)^2(n-m -1)} \sum_{1\le j < s \le n-m} (\xi^{(m)}_j -
\xi^{(m)}_s)^2.\eqno(17)$$

\begin{proof}  In fact,  the first Vi\'{e}te  formula (see \cite{Rah}) says that the coefficient $a_1$ ($a_0=1$) in (1) is equal to
$$- a_1= r_{1}\lambda_1  +  r_{2} \lambda_{2} +   \dots +  r_{k}\lambda_{k}.$$
On the other hand, differentiating (1) $n-1$ times, we find $z_{n-1}
= - a_1/ n$.  Thus in view of  (2) we prove the first identity  in (15). The second
identity  can be done similarly, using that a  centroid  is differentiation invariant, see, for instance, 
\cite{Rah}.  In order to establish the first identity  in (16), we
call formula (11) to find
$$\frac{f^{(n-2)} (z)}{(n-2)!} =  \frac{n(n-1)}{2} (z - z_{n-2}) (z + z_{n-2}- 2 z_{n-1}).\eqno(18)$$
Moreover,  as a consequence of the second Vi\'{e}te formula,  the
coefficient $a_2$ in (1), which  equals
$$a_2=   \frac{f^{(n-2)} ( z)}{(n-2)!} -  \frac{n(n-1)}{2} z^2 + n(n-1) z_{n-1}z\eqno(19)$$
can be expressed as follows
$$a_2=   {1\over 2} \left(\sum_{j=1}^k r_{j} \lambda_j\right)^2 -  {1\over 2} \sum_{j=1} ^k r_j\lambda_j^2.\eqno(20) $$
Hence letting $z=z_{n-2}$ in (18),  and taking into account (15) with $z=0$,
we deduce

$$2a_2=   n^2z^2_{n-1}  -   \sum_{j=1} ^k r_{j} \lambda_j^2 =  2 n(n-1) z_{n-1}z_{n-2}  -  n(n-1) z_{n-2}^2 . $$
Therefore,   using again (15) and (2), we easily come up  with the
first identity in (16). The second one can be proven  in the same
manner, involving roots of derivatives. Finally, we prove the first
identity in (17). Concerning the second identity,  see Lemma 6.1.5
in \cite{Rah}. Indeed, calling the first identity  in (16),  letting
$z= z_{n-1}$ and employing (15), we derive

$$  n^2(n-1)(z_{n-1} - z_{n-2})^2= n \sum_{j=1}^k r_{j}\lambda^2_j +
\left(\sum_{s=1}^k r_{s}\lambda_s\right)^2 - 2\left(\sum_{s=1}^k
r_{s}\lambda_s\right)\left(\sum_{j=1}^k r_{j}\lambda_j\right)$$
$$= n \sum_{j=1}^k r_{j}\lambda^2_j - \sum_{s=1}^k r^2_{j}\lambda^2_j -
2\sum_{1\le j <  s\le k} r_{j}r_{s}\lambda_j\lambda_s=  \sum_{1\le j
<  s\le k} r_{j}r_{s}(\lambda_j- \lambda_s)^2.$$

\end{proof}

The following result gives an identity, which is associated with  zeros of a monic
polynomial and common zeros  of its derivatives. Precisely, we have

{\bf Lemma 2.} {\it Let $f$ be a monic polynomial of exact degree $n\ge 2$, having $k$ distinct  roots of multiplicities $(2)$. Let  $z_{n-1}=\lambda_1$ be a common root of $f$ of multiplicity $r_1$ with  the unique root of its $n-1$st
derivative. Let also $z_m= \xi_{n-m}^{(m)}= \lambda_{k_m}$ be a common root of $f$ of multiplicity $r_{k_m}$ and
its $m$-th derivative, $m \in \{1,2,\dots,  n-2\}$. Then,  involving other roots of $f^{(m)}$,
the following identity holds}

$$\left[\frac{n-m-2}{(n-m)^2} + \frac{r_{k_m}+r_1-n}{n (n-1)}\right]\sum_{s=1}^{n-m-1} (z_m- \xi^{(m)}_s)^2
+ \frac{n-m-2}{(n-m)^2}\sum_{1\le s < t \le n-m-1} (\xi^{(m)}_s- \xi^{(m)}_t)^2
$$$$= \frac{(n-m)^2 r_{k_m}- (n-r_1)(n-m+2)}{n (n-1)} (z_m- z_{n-1})^2$$$$+ \
\frac{2}{n (n-1)}\sum_{j\neq 1, k_m} r_j \sum_{1 \le s < t \le
n-m-1} (\lambda_j- \xi^{(m)}_s)(\lambda_j- \xi^{(m)}_t).\eqno(21)$$

\begin{proof}  We begin, appealing to  (15) and letting  $z=0$.  We get
$$\sum_{s=1}^{n-m} \xi^{(m)}_s= (n-m) z_{n-1}, \quad   \xi_{n-m}^{(m)}=z_m.\eqno(22) $$
Hence via identities  (17) with $z=z_m$ we write the chain of equalities
$$\sum_{1\le s < t \le n-m} (\xi^{(m)}_s- \xi^{(m)}_t)^2=
\frac{(n-m-1)(n-m)^2}{n (n-1)} r_{k_m}(z_m- z_{n-1})^2 +
\frac{(n-m-1)(n-m)^2}{n (n-1)}\sum_{j\neq 1, k_m} r_j (\lambda_j -
z_{n-1})^2$$
$$=\frac{(n-m-1)(n-m)^2}{n (n-1)} r_{k_m}(z_m- z_{n-1})^2 +
\frac{n-m-1}{n (n-1)}\sum_{j\neq 1, k_m} r_j \left(\lambda_j - z_m+
\sum_{s=1}^{n-m-1} (\lambda_j- \xi^{(m)}_s) \right)^2$$
$$=\frac{(n-m-1)(n-m)^2}{n (n-1)} r_{k_m}(z_m- z_{n-1})^2 +
\frac{n-m-1}{n (n-1)}\left[\sum_{j\neq 1, k_m} r_j (\lambda_j -
z_m)^2 +  \sum_{j\neq 1, k_m} r_j \left(\sum_{s=1}^{n-m-1}
(\lambda_j- \xi^{(m)}_s) \right)^2\right.$$
$$\left.+ 2 \sum_{j\neq 1, k_m} r_j \sum_{s=1}^{n-m-1} (\lambda_j -
z_m)(\lambda_j- \xi^{(m)}_s)\right]= \frac{(n-m-1)(n-m)^2}{n (n-1)}
r_{k_m}(z_m- z_{n-1})^2$$
$$+ \frac{n-m-1}{n (n-1)}\left[(2(n-m)-1) \sum_{j\neq 1, k_m} r_j (\lambda_j -
z_m)^2 +  \sum_{j\neq 1, k_m} r_j \left(\sum_{s=1}^{n-m-1}
(\lambda_j- \xi^{(m)}_s) \right)^2\right.$$
$$\left.+ 2 \sum_{j\neq 1, k_m} r_j \sum_{s=1}^{n-m-1} (\lambda_j -
z_m)(z_m - \xi^{(m)}_s)\right]= \frac{(n-m-1)(n-m)^2}{n (n-1)}
r_{k_m}(z_m- z_{n-1})^2$$
$$+ \frac{n-m-1}{n (n-1)}\left[(2(n-m)-1) \sum_{j\neq 1, k_m} r_j (\lambda_j -
z_m)^2 +  \sum_{j\neq 1, k_m} r_j \left(\sum_{s=1}^{n-m-1}
(\lambda_j- \xi^{(m)}_s) \right)^2\right.$$
$$\left. - 2 (n-m)(n-r_1)(z_m-z_{n-1})^2 \right]= \frac{(n-m-1)}{n (n-1)}
\left((n-m)^2 r_{k_m}- n+r_1\right) (z_m- z_{n-1})^2$$
$$+ (n-m-1)(3(n-m)-2)(z_{n-1}- z_{n-2})^2  +  \frac{(n-m-1)(n-r_1)}{n (n-1)}\sum_{s=1}^{n-m-1}
(z_{n-1}- \xi^{(m)}_s)^2$$$$-    \frac{r_{k_m}(n-m-1)}{n
(n-1)}\sum_{s=1}^{n-m-1} (z_{m}- \xi^{(m)}_s)^2+ 2 \ \frac{n-m-1}{n
(n-1)}\sum_{j\neq 1, k_m} r_j \sum_{1 \le s < t \le n-m-1}
(\lambda_j- \xi^{(m)}_s)(\lambda_j- \xi^{(m)}_t).$$ Applying again
 (17), (22), we split the right-hand side of the latter
identity  in (17) in two parts, selecting the root $z_m$. Thus in the
same manner after straightforward calculations it becomes
$$\left[\frac{n-m-2}{(n-m)^2} + \frac{r_{k_m}+r_1-n}{n (n-1)}\right]\sum_{s=1}^{n-m-1} (z_m- \xi^{(m)}_s)^2
+ \frac{n-m-2}{(n-m)^2}\sum_{1\le s < t \le n-m-1} (\xi^{(m)}_s-
\xi^{(m)}_t)^2
$$$$= \frac{(n-m)^2 r_{k_m}- (n-r_1)(n-m+2)}{n (n-1)} (z_m- z_{n-1})^2+ \
\frac{2}{n (n-1)}\sum_{j\neq 1, k_m} r_j \sum_{1 \le s < t \le
n-m-1} (\lambda_j- \xi^{(m)}_s)(\lambda_j- \xi^{(m)}_t),$$
completing the proof of Lemma 2.

\end{proof}

{\bf Remark 1}.  It is easy to verify identity (21) for the least case $m=n-2$, when the double sums are
empty and $\xi^{(n-2)}_1= 2z_{n-1} -z_{n-2} $ (see above).

{\bf  Corollary 1.}  {\it A polynomial with only real roots of degree $n\ge 2$ is trivial, if and only if its  $n-2$nd  derivative has a double root}.

\begin{proof}  Indeed, necessity is obvious.  To prove sufficiency we see that since the $n-2$nd  derivative  has a double real root $x_{n-2}$,  it  is equal to the root $x_{n-1}$ of the $n-1$st  derivative. Therefore letting in (16)
$z= x_{n-1}$, we find that its left-hand side becomes zero and,
correspondingly, all squares in the right-hand side are zeros. This
gives the  conclusion that all roots are equal to $x_{n-1}$.

\end{proof}

{\bf  Corollary 2.}  {\it Let $f$ be an arbitrary  polynomial of
degree $n \ge 3$ with at least two distinct roots, whose $n-2$nd
derivative has a double root. Then it contains at least one  complex
root}.

\begin{proof} In fact, if all roots are real it is trivial via Corollary 1.

\end{proof}

 Evidently, each derivative up to $f^{(r-1)}$ of a  polynomial $f$ with only real roots,  where $r$ is the maximal  multiplicity,   shares a root with $f$. Moreover,  owing to the  Rolle  theorem all roots of $f^{(m)}, \  m=r,r+1,\dots, n-1$ are simple,  we have that a possible common root with $f$ is simple too (we note,  that a number of common roots does not exceed $k-2$, because  minimal and  maximal  roots cannot be zeros of $f^{(m)},\ m\ge r$).  This circumstance gives an immediate

 {\bf Corollary 3.} {\it There exists no non-trivial polynomial with only real roots, having two distinct zeros  and sharing a root  with at least one of its derivatives, whose order exceeds $r-1,\  r= \hbox{max}_{1\le j\le k}  (r_j)$. }

\begin{proof} Indeed, in the case of existence of such a polynomial,  these two distinct roots cannot be within zeros  of any derivative $f^{(m)},\  m > r$ owing to the Rolle theorem.  Moreover, if any of the two roots  are in common with roots of  $f^{(r)}$, its multiplicity is greater than $r$, which is impossible.

\end{proof}

We extend Corollary 3 to  three distinct real roots.  Precisely, it leads  to

 {\bf Corollary 4.} {\it There exists no non-trivial polynomial $f$ of degree $n \ge 3$ with only real roots, having three distinct zeros  and sharing a root  with its $n-2$nd  and $n-1$st  derivatives. }

\begin{proof}  Assume  such a polynomial exists and let's denote  its roots  $\lambda_1= x_{n-1}, \   \lambda_2= x_{n-2}$ and $\lambda_3$ of multiplicities  $r_1, \ r_2, \ r_3$, respectively. Hence employing identities (16), we write for this case
$$  (n^2- n- r_2) (x_{n-1} - x_{n-2})^2=  r_{3}(\lambda_3 - x_{n-1}) ^2.$$
In the meantime,  squaring  both sides of the first identity  in (15) for this case after simple modifications , we obtain
$$ r_2^2(x_{n-1} - x_{n-2})^2=  r_{3}^2 (\lambda_3 - x_{n-1}) ^2.$$
Hence, comparing with the  previous equality, we come out with the relation
$$ (n^2- n- r_2) r_3=  r_2^2.$$
But $n=r_1+r_2+r_3,\ r_j \ge 1, j=1,2,3.$  Consequently,
$$r_2^2 \ge  n(n-1)- r_2 >  (n-1)^2- r_2 \ge (r_1+r_2)^2-r_2 \ge r_2^2+ r_2+ r_1^2 > r_2^2,$$
which is impossible.
\end{proof}

{\bf Remark 2.}  If we omit the condition for $f$ to have a common root with the $n-2$nd  derivative in Corollary 4, it becomes false. In fact,  this circumstance can be shown  by the counterexample $f(x)= x^3-x.$

The following result deals with the case of 4 distinct roots.  We have,

{\bf Corollary 5.} {\it There exists no non-trivial polynomial $f$ of degree $n \ge 4$ with only real roots, having four  distinct zeros  and sharing a root  with its $n-2$nd  and $n-1$st  derivatives. }

\begin{proof}  Similarly to the previous corollary, we assume the existence of such a polynomial and call  its  roots
$\lambda_1= x_{n-1}, \   \lambda_2= x_{n-2}$ and $\lambda_3, \lambda_4$ of multiplicities $r_j, \  j=1,2,3,4$, respectively.
Hence the first  identity  in (16)  yields
$$  (n^2- n- r_2) (x_{n-1} - x_{n-2})^2=  r_{3}(\lambda_3 - x_{n-1}) ^2+   r_{4}(\lambda_4 - x_{n-1}) ^2 .\eqno(23)$$
Meanwhile, using  the first identity  in (15) for this case, we derive in a similar manner
$$ r_2^2(x_{n-1} - x_{n-2})^2=  r_{3}^2 (\lambda_3 - x_{n-1}) ^2+ r_{4}^2 (\lambda_4 - x_{n-1}) ^2+
2 r_{3}r_4 (\lambda_3 - x_{n-1}) (\lambda_4 - x_{n-1}).$$
Thus,  after straightforward calculations,  we come out with the quadratic equation
$$Ay^2+ By+ C=0$$
in the variable $y=  (\lambda_3 - x_{n-1}) /  (\lambda_4 - x_{n-1})$ with coefficients $A= r_3r_2^2- r_3^2(n^2-n-r_2),\
B= - 2 r_3r_4 (n^2-n-r_2),\ C= r_4r_2^2- r_4^2(n^2-n-r_2).$   But, it is easy to verify that $B^2-4AC >0.$ Therefore the quadratic
equation has two distinct real roots.  Writing $\lambda_3 - x_{n-1}= y (\lambda_4 - x_{n-1})$ and substituting into (23), we obtain
$$ (n^2- n- r_2) (x_{n-1} - x_{n-2})^2= ( r_{3}y^2+ r_4) (\lambda_4 - x_{n-1})^2.$$
At the same time, since $y\neq 0$,  we have $  \lambda_4 - x_{n-1}=
y^{-1}  (\lambda_3 - x_{n-1})$ and
$$y^2 (n^2- n- r_2) (x_{n-1} - x_{n-2})^2= ( r_{3}y^2+ r_4) (\lambda_3 - x_{n-1})^2.$$
Hence,
$$\lambda_4= x_{n-1}  \pm \sqrt{\frac{n^2- n- r_2}{ r_{3}y^2+ r_4}} \  |x_{n-1} - x_{n-2}|,$$
$$\lambda_3= x_{n-1}  \pm |y| \sqrt{\frac{n^2- n- r_2}{ r_{3}y^2+ r_4}}\   |x_{n-1} - x_{n-2}|.$$
Consequently,
$$\lambda_4- \lambda_3 =  \sqrt{\frac{n^2- n- r_2}{ r_{3}y^2+ r_4}} \  |x_{n-1} - x_{n-2}| ( 1-|y|)=
 - \sqrt{\frac{n^2- n- r_2}{ r_{3}y^2+ r_4}} \  |x_{n-1} - x_{n-2}| ( 1+|y|)$$
 $$=  \sqrt{\frac{n^2- n- r_2}{ r_{3}y^2+ r_4}} \  |x_{n-1} - x_{n-2}| ( 1+|y|)=
 \sqrt{\frac{n^2- n- r_2}{ r_{3}y^2+ r_4}} \  |x_{n-1} - x_{n-2}| (|y|-1),$$
which is  possible only in the case $x_{n-1}=x_{n-2},$\  $\lambda_3=\lambda_4$. Thus we get a contradiction with Corollary 1 and complete the proof.
\end{proof}

In the same manner we prove

{\bf Corollary 6.} {\it There exists no non-trivial polynomial $f$ of degree $n \ge 5$ with only real roots, having five   distinct zeros  and sharing  roots  with its $n-2$nd  and $n-1$st  derivatives. }

\begin{proof}  Assuming its existence, it has  the roots  $\lambda_1= x_{n-1}, \   \lambda_2= x_{n-2}$,  $\lambda_3= 2x_{n-1}- x_{n-2}, \  \lambda_4$ and $\lambda_5$ of multiplicities $r_j, \  j=1,2,3,4, 5$, respectively. Hence
$$  (n^2- n- r_2-r_3) (x_{n-1} - x_{n-2})^2=  r_{4}(\lambda_4 - x_{n-1}) ^2+   r_{5}(\lambda_5 - x_{n-1}) ^2 .$$
Therefore  using similar ideas as in the proof of Corollary 5, we come out again to the contradiction.

\end{proof}

For an arbitrary  number of distinct zeros we establish the following

{\bf Corollary 7.} {\it There exists no non-trivial polynomial $f$ of degree $n$ with only real roots, having $k \ge 2$   distinct zeros of multiplicities $(2)$ $r_j,\  j=1,\dots, k$ and  among them all roots of $f^{(m)}$ for some  $m$, satisfying   the relations
$$ r \le m <  {1\over 2}\left(1-{1\over r_0}\right)(n-1),\eqno(24)$$
where $r,\  r_0$ are maximum and minimum multiplicities of roots of $f$.}

\begin{proof}  In fact, as a consequence of (16) we have the identity 
$${(n-m)(n-m -1)\over n(n-1)}  \sum_{j=1}^k r_{j}(\lambda_j- x_{n-1})^2 =  
 \sum_{j=1}^{n-m} (\xi^{(m)}_j - x_{n-1})^2\eqno(25)$$
for some $m$,  satisfying condition (24).  Hence, since $m\ge r$,  it has $n-m \le k-2$ and $\xi^{(m)}_j= \lambda_{m_j}, \  m_j \in \{1,\dots, k\}, \  j= 1,\dots, n-m$  are simple roots of $f^{(m)}$.  Thus  we find 
$$ \sum_{j=1}^{n-m} \left[ r_{m_j} {(n-m)(n-m -1)\over n(n-1)} -1\right]  (\lambda_{m_j}- x_{n-1})^2 + 
{(n-m)(n-m -1)\over n(n-1)}   \sum_{j=n-m+1}^k r_{m_j}(\lambda_{m_j} - x_{n-1})^2 =  0.$$
But, owing to condition (24) 
$$r_{m_j} {(n-m)(n-m -1)\over n(n-1)} -1 \ge  r_{0} {(n-m)(n-m -1)\over n(n-1)} -1 \ge 0, \ j=  1,\dots, n-m.$$
Indeed, we have from the latter inequality 
$$m \le n- {1\over 2} - \sqrt{\frac{n^2-n}{r_0} + {1\over 4}}$$
and, in turn,
$$ n- {1\over 2} - \sqrt{\frac{n^2-n}{r_0} + {1\over 4}}=  \frac{2 (1-r_0^{-1}) (n^2-n)}
{ 2n- 1  + \sqrt{4(n^2-n)r^{-1}_0 + 1}}\ge \frac{ (1-r_0^{-1}) (n^2-n)}
{ 2n- 1} >  {1\over 2}\left(1-{1\over r_0}\right)(n-1).$$
Therefore  $\lambda_j= x_{n-1}, \ j=1,\dots, k$ and this contradicts to the fact that all roots are distinct.  
\end{proof}

Finally, in this section, we will employ identities (17) to prove an analog of the Obreshkov- Chebotarev theorem for multiple roots (see \cite{Rah}, Theorem 6.4.3), involving estimates for smallest and largest of distances between consecutive zeros of polynomials and their derivatives.  Namely, it has

{\bf Theorem 2.} {\it  Let $f$ be a polynomial of degree $n >  2$ with only real zeros. Denote the largest and the smallest of the distances between consecutive zeros of $f$ by $\Delta$ and $\delta$, respectively.  Denoting the corresponding quantities associated with $f^{(m)}, \ m=1,2,\dots,\  n-2$ by  $\Delta^{(m)}$ and $\delta^{(m)}$, the following inequalities hold

$$\delta^{(m)}  \le \Delta \  {rk \over n} \  \sqrt{  \frac{k^2-1}{ (n-m+1)(n-1)}},\eqno(26)$$
$$\delta \  {r_0 k \over n} \  \sqrt{  \frac{k^2-1}{ (n-m+1)(n-1)}}\le  \Delta^{(m)} ,\eqno(27)$$
$$\delta \  {r_0 k \over 2 n} \  \sqrt{  \frac{k^2-1}{ 3 (n-1)}}\le |x_{n-1} - x_{n-2}| \le
 \Delta \  {rk \over 2 n} \  \sqrt{  \frac{k^2-1}{ 3(n-1)}},\eqno(28)$$
where $r_0,\  r$ are minimum and maximum multiplicities of roots of $f$, respectively,  and $k \ge 2$ is a number of distinct roots.}

\begin{proof}  Following similar ideas as in the proof of Theorem 6.4.3 in \cite{Rah},  we assume distinct roots of $f$ in the increasing order and roots of its  $m$-th derivative in the non-decreasing order, and taking the second identity in (17), we deduce
$$ {[\delta^{(m)}]^2 \over (n-m)^2(n-m -1)} \sum_{1\le j < s \le n-m} (s-j)^2 \le {[\Delta r]^2  \over n^2(n-1)} \sum_{1\le j < s \le k} (s-j)^2.$$
Hence, in view of  the value of the sum
$$ \sum_{1\le j < s \le q} (s-t)^2 = {1\over 12} q^2(q^2-1),$$
after simple manipulations we arrive at the inequality (26).  In the same manner (cf. \cite{Rah}) we establish inequalities (27), (28), employing Sz.-Nagy type identities (17). 
\end{proof}

\section{Laguerre  type inequalities }

In 1880 Laguerre proved his famous theorem for polynomials with only
real roots, which provides their localization with  upper and lower
bounds (see details in \cite{Rah}).  Precisely,  we have the following
Laguerre inequalities
$$x_{n-1}-  (n-1)\left|x_{n-1}- x_{n-2}\right| \le w_j \le  x_{n-1} + (n-1)\left|x_{n-1}- x_{n-2}\right|, \
 j=1,\dots, n,$$
where $w_j$ are  roots of the polynomial $f$ of degree $n$ and $x_{n-1}, \  x_{n-2}$ are roots of  $f^{(n-1)},\ f^{(n-2)}$, respectively.   First  we prove an analog of the Laguerre inequalities  for multiple roots.

{\bf  Lemma 3.}  {\it Let $f$ be a polynomial with only real roots of degree  $n \in \mathbb{N}$,  having $k$  distinct roots
 $\lambda_j, \ j=1,\ \dots, k$ of multiplicities $(2)$ and $x_{n-1}, \  x_{n-2}$ be roots of  $f^{(n-1)},\ f^{(n-2)}$, respectively.
  Then the following Laguerre   type inequalities  hold}
$$x_{n-1}-  \sqrt{\frac{(n-r_j)(n-m-1)}{r_j-m} }\left|x_{n-1}- x_{n-2}\right| \le \lambda_j \le  x_{n-1} +
\sqrt{\frac{(n-r_j)(n-m-1)}{r_j-m} }\left|x_{n-1}- x_{n-2}\right|,\eqno(29)$$
where $ j=1,\dots, k, \  m= 0,1,\dots, r_j-1.$

\begin{proof} In fact, appealing to the  Sz.-Nagy type identities (15), (16) and the Cauchy -Schwarz inequality, we find
$$   (x_{n-1} - x_{n-2})^2={1\over (n-m)(n-m-1)} \left[ \sum_{s=1}^{n-m} (\xi^{(m)}_s-  \lambda_j)^2-  (n-m) (x_{n-1} -
\lambda_j)^2\right]$$$$ \ge  {1\over (n-m)(n-m-1)} \left[
\frac{1}{n-r_j} \left( \sum_{s=1}^{n-m} (\xi^{(m)}_s-
\lambda_j)\right)^2-  (n-m) (x_{n-1} - \lambda_j)^2\right]$$$$=
\frac{r_j-m} {(n-r_j)(n-m-1)} \left(x_{n-1}- \lambda_j\right)^2, \
m= 0,1,\dots, r_j-1,$$ which yields (29).

\end{proof}

As a corollary we improve the Laguerre inequality (28) for multiple roots.

{\bf Corollary 8.} {\it  Let $f$ be a polynomial with only real roots of degree  $n \in \mathbb{N}$.
Then the multiple zero $\lambda_j$ of multiplicity $r_j\ge 1, \   j=1,\dots, k$ lies in the interval}
$$\left[ x_{n-1}-  \sqrt{\left(\frac{n}{r_j}-1\right)(n-1) }\left|x_{n-1}- x_{n-2}\right|, \quad   x_{n-1} +
\sqrt{\left(\frac{n}{r_j}-1\right)(n-1) }\left|x_{n-1}- x_{n-2}\right| \right].\eqno(30)$$

\begin{proof}  Indeed, the fraction $\frac{(n-r_j)(n-m-1)}{r_j-m}$ attains its minimum value, letting $m=0$ in (29).
\end{proof}

{\bf Remark 3.}  When all roots are simple, the latter interval
coincides with the one generated by (28).

A localization of roots of the  $m$-th derivative $f^{(m)}, \ m=0,1,\dots,  n-2$ is given by

{\bf  Lemma 4.}  {\it Roots of the $m$-th derivative $f^{(m)}, \ m=0,1,\dots,  n-2$ satisfy the following Laguerre   type inequalities}
$$x_{n-1}-  (n-m-1)\left| x_{n-1}- x_{n-2}\right| \le \xi^{(m)}_\nu \le
 x_{n-1} +  (n-m-1)\left| x_{n-1}- x_{n-2}\right|,\eqno(31)$$
where   $\nu=1,\dots, n-m.$

\begin{proof} Similarly to the proof of Lemma 3, we employ  the  Sz.-Nagy type identities (15), (16) and the Cauchy -Schwarz inequality to deduce
$$   (x_{n-1} - x_{n-2})^2={1\over (n-m)(n-m-1)} \left[ \sum_{s=1}^{n-m} (\xi^{(m)}_s-  \xi^{(m)}_\nu)^2-  (n-m) (x_{n-1} -
\xi^{(m)}_\nu)^2\right]$$$$ \ge  {1\over (n-m)(n-m-1)} \left[
\frac{1}{n-m-1} \left( \sum_{s=1}^{n-m} (\xi^{(m)}_s-
\xi^{(m)}_\nu)\right)^2-  (n-m) (x_{n-1} -
\xi^{(m)}_\nu)^2\right]$$$$= \frac{1}{(n-m-1)^2} \left(x_{n-1}-
\xi^{(m)}_\nu\right)^2, \  m= 0,1,\dots, n-2.$$ Thus we come up
with (31) and complete the proof.

\end{proof}

When the root $x_{n-1}=\lambda_1$ be in common with $f$ of multiplicity $r_1$, we have

{\bf  Lemma 5.}  {\it Let $f$ be a polynomial with only real roots
of degree $n \ge 2$ and $x_{n-1}=\lambda_1$ be a  common zero with
$f$ of multiplicity $r_1$, having $k \ge 2$  distinct roots $\lambda_j$ of multiplicities
$r_j, j=1,\dots, k$. Then the following Laguerre type inequalities
hold}
$$x_{n-1}-   \sqrt{\left({1\over r_s}- {1\over n-r_1}\right) (n^2-n)}\left| x_{n-1}- x_{n-2}\right| \le \lambda_s \le  x_{n-1}
 $$$$+  \sqrt{\left({1\over r_s}- {1\over n-r_1}\right) (n^2-n)}\left| x_{n-1}- x_{n-2}\right|,\eqno(32)$$
where   $s=2,\dots, k.$

\begin{proof} In the same manner  we involve   the  first Sz.-Nagy type identity in  (15) with $z=  \lambda_s$,  which can be written in the form
$$(n-r_1) (x_{n-1} - \lambda_s)  =  \sum_{j=2}^k r_{j}(\lambda_j- \lambda_s).$$
Hence squaring  both sides of the latter equality and appealing to the Cauchy -Schwarz inequality, we derive by virtue of  (16)
$$(n-r_1)^2 (x_{n-1} - \lambda_s)^2  = \left( \sum_{j=2}^k r_{j}(\lambda_j- \lambda_s)\right)^2$$
$$\le (n-r_1-r_s)\sum_{j=2}^k r_{j}(\lambda_j- \lambda_s)^2=  (n-r_1-r_s)\left[ (n^2-n) ( x_{n-1}- x_{n-2})^2 +
(n-r_1)(x_{n-1}- \lambda_s)^2\right].$$
Thus after simple calculations we easily arrive at (32).

\end{proof}

{\bf Remark 4}. Inequalities (27) are sharper than the corresponding relations, generated by interval (30).

The following result gives a Laguerre type localization for common roots of a possible CA-polynomial with only real roots and its $m$-th derivative.

{\bf Lemma 6.} {\it     Let $f$ be a CA-polynomial of degree $n \ge
2$ with only real distinct zeros of multiplicities $(2)$,   including common
roots $x_{n-1}=\lambda_1$ of its $n-1$st   derivative and $x_m$ of
its $m$-th derivative,   $m=  r, r+1, \dots, n-2$, where  $r=
\hbox{max}_{1\le j\le k} (r_j)$. Then the following Laguerre type
inequality holds
$$ \frac{n-r_1- r_{k_m}}{(n-r_1)^2}
\left(n^2-r_1+ (n-r_1)(n-m) (n-m-2)\right) (x_{n-1}-  x_{n-2})^2 \ge  (x_{n-1}- x_{m})^2,\eqno(33)$$
where $x_{n-2}$ is a root of $f^{(n-2)}$ and $r_{k_m}$ is the multiplicity of  $x_m$ as a root of $f$}.

\begin{proof}   Appealing  again to  Sz.-Nagy  type  identities (15), (16) with $z=x_m$, inequality (31)  and the Cauchy-Schwarz inequality, we  find
$$   (x_{n-1} - x_{n-2})^2={1\over n(n-1)} \left[ \sum_{j=2}^k r_{j}(\lambda_j- x_m)^2-  (n-r_1) (x_{n-1} -
x_m)^2\right]$$$$ \ge {1\over n(n-1)} \left[ \sum_{j=2}^k
r_{j}(\lambda_j- x_m)^2- (n-r_1)(n-m-1)^2(x_{n-1}-  x_{n-2})^2\right]$$
$$ \ge {1\over n(n-1)} \left[ {1\over n-r_1- r_{j_m}} \left(\sum_{j=2}^k
r_{j}(\lambda_j- x_m)\right)^2 -  (n-r_1)(n-m-1)^2(x_{n-1}-  x_{n-2})^2\right]$$
$$= {n-r_1\over n(n-1)} \left[ {n-r_1\over n-r_1- r_{j_m}} (x_{n-1}- x_m)^2 - (n-m-1)^2(x_{n-1}-  x_{n-2})^2\right].$$
Hence, making straightforward calculations,  we  derive (33), completing the proof of Lemma 6.

\end{proof}

Let us denote by $d,\ d^{(m)},\  D, D^{(m)}$ the following values
$$   d=  \hbox{min}_{2\le j\le k}  |\lambda_j-   x_{n-1}|, \quad
  d^{(m)}=  \hbox{min}_{1\le j\le n-m}   |\xi^{(m)}_j -   x_{n-1}|,\eqno(34)$$
$$   D=  \hbox{max}_{2\le j\le k}  |\lambda_j-   x_{n-1}|, \quad
  D^{(m)}=  \hbox{max}_{1\le j\le n-m}   |\xi^{(m)}_j -   x_{n-1}|,\eqno(35)$$
and by
$$\hbox{span} (f) = \lambda^*- \lambda_*,$$
where
$$  \lambda^*= \hbox{max}_{1\le j\le k} (\lambda_j), \quad     \lambda_*= \hbox{min}_{1\le j\le k} (\lambda_j)$$
are roots of $f$ with  multiplicities $r^*, \ r_*$, respectively.   Then   $ D^{(m+1)}  \le  D^{(m)}\le D$ and
  (cf. \cite{Rah})  $\hbox{span}(f^{(m+1)}) \le \hbox{span}(f^{(m)}) \le \hbox{span}(f)$, where $\hbox{span}(f^{(m)})$
   is the span of the $m$-th derivative.  Moreover, the strict inequalities $ D^{(m)}  <  D$,\    $\hbox{span}(f^{(m)}) < \hbox{span}(f)$   hold when $m$ is sufficiently large.

  {\bf Lemma 7}. {\it Let $x_{n-1}=\lambda_1, \  x_{n-2}=\lambda_2$ be  common roots of $f$  with its $n-1$st, $n-2$nd  derivatives, respectively,  of multiplicities $r_1, r_2$ as roots of $f$,  and the maximum distance $D$ (see $(35)$) be attained at the root $\lambda_{s_0}, \ s_0 \in \{ 3,\dots, k\},\ k \ge 3$   of $f$ of multiplicity $r_{s_0}$.   Then the following inequalities hold }
    $$\sqrt{\frac{n^2-n- r_2}{n-r_1-r_2}}\  \left|x_{n-1} - x_{n-2}\right|   \le D \le \sqrt{\frac{n^2-n-r_2}{r_{s_0}}}
     \left|x_{n-1} - x_{n-2}\right|,\eqno(36)$$
  $${1\over 2} \sqrt{ {r_{s_0} \over 3(n-r_1) }\left(5 + \frac{ r_2} {n^2-n-r_2}\right)}\hbox{span}(f) \le D \le \sqrt{{1\over n-r_1} \left[ n-r_1- {r_{s_0} \over 4}\left(5 + \frac{ r_2} {n^2-n-r_2}\right)\right] } \hbox{span}(f).\eqno(37)$$
 
\begin{proof}  In order to establish (36), we employ identities (16) and under condition of the lemma we write 
$$ (n^2-n- r_2)(x_{n-1} -   x_{n-2})^2 = \sum_{j=3}^k r_{j}(\lambda_j- x_{n-1})^2 \le (n-r_1-r_2) D^2.$$
Since $n > r_1+r_2$ and $x_{n-2}\neq \lambda_{s_0}$ (otherwise $f$ is trivial, because equalities $x_{n-2}= \lambda_{s_0}= \lambda^*$ or  $x_{n-2}= \lambda_{s_0}= \lambda_*$ mean that the maximum multiplicity $r > n-2$,  and we appeal to Corollary 3),  we come up with the lower  bound (36) for $D$.  The lower bound comes immediately from the estimate 
$$(n^2-n- r_2)(x_{n-1} -   x_{n-2})^2 = \sum_{j=3}^k r_{j}(\lambda_j- x_{n-1})^2    \ge  r_{s_0} D^2.$$
Now,  since $2D \ge  \hbox{span}(f)$, we find from (36)
$$\hbox{span}(f) \le 2 \sqrt{\frac{n^2-n-r_2}{r_{s_0}}} \  \left|x_{n-1} - x_{n-2}\right|.$$
Hence,  since  $D = \hbox{max} \left(|\lambda^*- x_{n-1}|, \ |\lambda_*- x_{n-1}|\right)$, the
$n-2$nd  derivative has roots $x_{n-2}$ and $2x_{n-1}- x_{n-2}$ and
$\hbox{span}(f)= D + \Lambda$, where $\Lambda = \hbox{min}
\left(|\lambda^*- x_{n-1}|, \ |\lambda_*- x_{n-1}|\right)$, we
appeal to the first identity  in (16), letting $z= \lambda_{s_0}$ and writing it in the form
$$ (n-r_1)(x_{n-1} -  \lambda_{s_0})^2 = \sum_{j=2}^k r_{j}(\lambda_j- \lambda_{s_0})^2 - n(n-1)(x_{n-1} -
x_{n-2})^2.$$
Therefore,
$$(n-r_1)D^2 \le \left[ n-r_1- {5\over 4} r_{s_0} - \frac{r_{s_0} r_2} {4(n^2-n-r_2)}\right]  [\hbox{span}(f)]^2$$ 
and we establish the upper bound (37) for $D$. On the other hand $\hbox{span}(f)= D+ \Lambda$. So,
$$D^2 \le \left(1 -   {r_{s_0} \over 4(n-r_1) }\left(5 + \frac{ r_2} {n^2-n-r_2}\right)\ \right) \left(D^2 + \Lambda^2 +
2D\Lambda\right)$$ and we easily come out with the lower bound (37)
for $D$, completing the proof of Lemma 7.

\end{proof}

{\bf Lemma 8}. {\it Let $x_{n-1}=\lambda_1, x_{n-2}=\lambda_2$ be
common roots of $f$ with its $n-1$st, $n-2$nd    derivatives of
multiplicities $r_1, r_2,\ r_1+r_2 < n$, respectively. Then  we
 have the following lower bound for $\hbox{span}(f)$}
  $$\hbox{span}(f)\ge \sqrt{\frac{n^2-r_1}{n-r_1-r_2}}\ |x_{n-1}-x_{n-2}|.\eqno(38)$$

\begin{proof} Indeed,  identities  (16) with $z=x_{n-2}$ yield
$$(n^2-r_1)(x_{n-1}-x_{n-2})^2=\sum_{j=3}^k r_{j}(\lambda_j- x_{n-2})^2$$
and we derive
$$(n^2-r_1)(x_{n-1}-x_{n-2})^2\le (n-r_1-r_2)[\hbox{span}(f)]^2,$$
which implies (38).
\end{proof}

Next, we establish an analog of Lemma 5 for roots of derivatives.   Precisely,  we have

 {\bf Lemma 9}. {\it Let  $x_{n-1},  \  x_{n-2}$ be roots of the $n-1$st, $n-2$nd    derivatives of $f$, respectively.  Then
 $$D^{(m)} \ge \sqrt{n-m-1}\  |x_{n-1}-x_{n-2}|,\eqno(39)$$
 where $m \in \{r, r+1,\dots, n-2 \}, \   r= \hbox{max}_{1\le j\le k} (r_j).$  Besides, if $x_{n-1}$ is a root of $f^{(m)}$, then
 we have a stronger inequality $$D^{(m)} \ge \sqrt{n-m}\  |x_{n-1}-x_{n-2}|.\eqno(40)$$
Moreover,
  $$2\  D^{(m)} \ge \hbox{span}(f^{(m)}) \ge   \frac{n-m}{n-m-1}\   D^{(m)}.\eqno(41)$$
  and if $x_{n-1}$ is a  root of  $f^{(m)}$,  it becomes
   $$2\  D^{(m)} \ge \hbox{span}(f^{(m)}) \ge  \sqrt{ \frac{(n-m)(n-m-1)+1}{(n-m-1)(n-m-2)}}\   D^{(m)},\eqno(42)$$
where  $m \in \{r, r+1,\dots, n-3 \}.$}
\begin{proof}  In fact,  since (see (16))
$$(n-m)(n-m-1)(x_{n-1} - x_{n-2})^2= \sum_{j=1}^{n-m} (\xi^{(m)}_{j} -  x_{n-1})^2 \le (n-m) \left[ D^{(m)}\right]^2,$$
we get (39). Analogously, we immediately come out with (40), when $x_{n-1}$ is a root of $f^{(m)}$, because one element of the sum of squares is zero.  In order to prove (41), we appeal again to (16), letting  $z= \xi^{(m)}_{s_0}, \ s_0 \in \{1,2, \dots,  n-m\}$,  $m \in \{r, r+1,\dots, n-2 \}, \   r= \hbox{max}_{1\le j\le k} (r_j)$, which is a root of  $f^{(m)}$, where the maximum $D^{(m)}$ is attained.  Hence owing to Laguerre type inequality (31)
$$(n-m)\left[D^{(m)}\right]^2\le (n- m-1) [\hbox{span}(f^{(m)})]^2 - \frac{n-m}{n-m-1} \left[D^{(m)}\right]^2,$$
which leads  to  the lower bound for $ \hbox{span}(f^{(m)})$ in   (41).   The upper bound is straightforward
since $x_{n-1}$ belongs to the smallest interval containing roots of $f^{(m)}$.  In the same manner we establish (42),  since in this case
$$(n-m-1)\left[D^{(m)}\right]^2\le (n- m-2) [\hbox{span}(f^{(m)})]^2 - \frac{n-m}{n-m-1} \left[D^{(m)}\right]^2.$$

\end{proof}

{\bf Remark 5}.  The case $m=n-2$ gives  equalities in (39), (41).  Letting the same value of $m$ in (40), we easily get
a contradiction, which means that the only trivial polynomial is within polynomials with only real roots,
whose derivatives $f^{(n-2)}, \  f^{(n-1)}$ have a common root (see Corollary 1).

\section{Applications to the Casas- Alvero conjecture}

In this final section we will discuss properties of possible CA-polynomials,  which share roots with each of their non-constant derivatives. We will investigate particular cases of the Casas-Alvero conjecture, especially for polynomials with only real roots, showing when it holds true or, possibly, is false. 

We begin with

{\bf Proposition 1}. {\it The Casas-Alvero conjecture holds true, if and only if  it is true for common roots $\{z_\nu \}_0^{n-1}$ lying in the unit circle.}

\begin{proof}   The necessity is trivial. Let's   prove the sufficiency. Let the conjecture be true for common roots $\{z_\nu \}_0^{n-1}$ of a complex polynomial $f$ and its non-constant derivatives, which lie in the unit circle.  Associating with $f$ an Abel-Goncharov polynomial $G_n$ (6),    one can choose an arbitrary $\alpha  >0$ such that $\ |z_\nu| < \alpha^{-1}, \   \nu=0, 1, \dots , n-1. $  Hence owing to  (7)
$$ f (\alpha z_\nu) = G_n\left(\alpha z_0, \alpha z_\nu, \alpha z_1,\dots, \alpha z_{n-1}\right)
= \alpha^n G_n(  z_\nu) = \alpha^n f (  z_\nu)=  0,\  \nu=0, 1,\dots , n-1,$$
and
$$f^{(\nu)} _n(\alpha z ) =  n! {d^{\nu}\over d z^{\nu}}
\int_{\alpha z_0}^{\alpha z}  \int_{\alpha z_1}^{s_{1}} \dots \int_{\alpha z_{n-1}}^{s_{n-1}} d s_{n} \dots d s_{1}
=n!  \alpha ^\nu \int_{\alpha z_\nu}^{\alpha z}  \int_{\alpha z_{\nu+1}}^{s_{\nu+ 1}} \dots \int_{\alpha z_{n-1}}^{s_{n-1}} d s_{n} \dots d s_{\nu+1},$$
we find $f^{(\nu)} _n(\alpha z_\nu ) =0$.   Hence  $\alpha z_\nu, \  \nu=0, 1,\dots , n-1$ are common
roots of  $\nu$-th derivatives  $f^{(\nu)}$ and $f$,  lying in the unit circle.   Consequently,  since via assumption the
Casas-Alvero conjecture is true when common roots are inside the unit circle, we have that  $f$ is trivial  and $z_0=z_1=\dots = z_{n-1} = a$ is a unique joint root of $f$ of the multiplicity $n$.  Proposition 1  is proved.
\end{proof}

The following lemma will be useful in the sequel.

{\bf Lemma 10}. {\it Let  $f$ be a  CA-polynomial with only real roots of degree $n \ge 2$ and  $\{x_\nu
\}_{0}^{n-1}$ be a sequence of common  roots of $f$ and the corresponding derivatives $f^{(\nu)}$.  
Let  $f^{(s+\nu)} (x_\nu) \ge 0, \  s =1,2,\dots, n-\nu-1$ and $\nu=0,1,\dots, n-1$. Then $x_\nu$ is a
maximal root of the derivative  $f^{(\nu)}$.}

\begin{proof}  In fact,  the proof is an immediate consequence of  the expansion (12), where we let $G_n(x)=f(x)$.  Indeed,  $f^{(\nu)} (x_\nu) =0,   \nu=0,1,\dots, n-1$ and when $ x > x_\nu$   we have from (12) $f^{(\nu)} (x) > 0,   \nu=0,1,\dots, n-1$.   So,  this  means that there are  no roots, which are bigger than $x_\nu$. This completes  the proof of Lemma 10.
\end{proof}

{\bf Proposition  2}. {\it  Under conditions of Lemma 10  the Casas-Alvero conjecture holds true for polynomials with only real roots.}

\begin{proof}  We will show that under conditions of Lemma 10 there exists no CA-polynomial $f$ with only real roots.  Indeed, assuming its existence,  we find via conditions of  the lemma that the root $x_0$ is a maximal zero  of $f(x)$. This means that $x_0 \ge x_1$. On the other hand, the classical  theorem of Rolle states that between zeros $x_0, \  x_1$ in the case $x_0 > x_1$ there exists at least one zero of the derivative $f^{(1)} (x)$, say $\xi_1^{(1)}$,
 which is bigger than $x_1$.  But this is impossible because $x_1$ is a maximal zero of the first derivative.
  Thus $x_0=x_1\ge x_2$.   Then   between $x_1$ and $x_2$ in the case $x_1 > x_2$ there exists a zero $\xi_2^{(1)}$
   of the first derivative such that $x_1>  \xi_2^{(1)}  > x_2$.  Hence between $x_1$ and   $\xi_2^{(1)}$ there exists at least one     zero of the second derivative, which is bigger than $x_2$. But this is impossible, since $x_2$ is a maximal zero of $f^{(2)} (x)$.      Therefore   $x_0=x_1=x_2$.   Continuing this process we observe that the sequence  $\{x_\nu \}_0^{n-1}$ is stationary and $f$ has a unique joint root, which contradicts the definition of the CA-polynomial. 
\end{proof}

{\bf Corollary 9}. {\it   There exists no CA-polynomial $f$ with only real roots, having a non-increasing sequence $\{x_\nu \}_0^{n-1}$ of roots in common with $f$ and its non-constant derivatives.}

\begin{proof}    Obviously,  via (13) $f^{(s+\nu)} (x_\nu) \ge 0, \  s=1,2,\dots, n-\nu-1$ and conditions of Lemma 10 are satisfied.
\end{proof}

{\bf Corollary 10}. {\it   There exists no CA-polynomial $f$ with only real roots, such that each $x_\nu$ in the  sequence $\{x_\nu \}_0^{n-1}$  is a maximal root of the derivative  $f^{(\nu)} (x),  \  \nu=0,1,\dots, n-1$.}

\begin{proof}    The proof is similar to the proof of Proposition  2.
\end{proof}

An immediate consequence of Corollaries 3,4,5 is 

{\bf Corollary 11}. {\it   The CA-polynomial, if any,   with only real roots has at least 5 distinct zeros. }

 Let us  denote by $l(m)$ the  number of distinct roots of the $m$-th derivative $f^{(m)},\ m=0, 1,\dots, n-2$,   which are in common with $f$ and different from $\lambda_1=x_{n-1}$, which is a common root with $f^{(n-1)}$,   i.e.  the $m$-th derivative $f^{(m)}$ has $l(m)$ common roots with $f$

$$\lambda_{j_1}, \dots,  \lambda_{j_{l(m)}}  \subseteq  \{ \lambda_2, \lambda_3, \dots, \  \lambda_k\} $$
of multiplicities
$$r_{j_1}, \dots,  r_{j_{l(m)}}  \subseteq  \{ r_2, r_3, \dots, \  r_k\} .$$
For instance,  $l(0)= k-1, \   l(1)= k-1-s$, where $s$ is a number
of simple roots of $f$. So, we see that $n-m \ge l(m) \ge 0$ and
since $f$ is a  CA-polynomial,  $l(m)=0$ if and only if
$x_{n-1}=\lambda_1$ is the only common root of $f$ with $f^{(m)}$.

{\bf Lemma  11}. {\it There exists no  CA-polynomial with only real roots, having the property
$l(m)= l(m+1)=0$ for some $m \in \{r, r+1, \dots,n-2\},$ where $r= \hbox{max}_{1\le j\le k} (r_j)$.}

\begin{proof} In fact, as we saw above, since all roots are real,  it follows that all roots of $f^{(m)}, \  m \ge r$ are simple, which contradicts equalities $l(m)= l(m+1)=0$.   Indeed, the latter equalities yield that $x_{n-1}$ is a multiple root of $f^{(m)}$. Therefore $r\ge r_1 >  m+1\ge r+1$, which is impossible.
\end{proof}

Further, as in Lemma 7 we involve the root  $\lambda_{s_0}$ of multiplicity $r_{s_0}$,  and $D= |\lambda_{s_0}- x_{n-1}|$ (see (35)). Thus $\lambda_{s_0}= \lambda_*$ or $\lambda_{s_0}= \lambda^*$ and, correspondingly,  $r_{s_0}=r_*$ or $r_{s_0}=r^*$. Hence, calling  Sz.-Nagy identities (15), we  let  $z= x_{n-1}$ and assume without loss of generality that $\lambda_{s_0}= \lambda^*$. Then  we obtain for  $m \ge r$
$$r_*(x_{n-1}- \lambda_*) = r^*D + \sum_{j=2, \ r_j\neq r_*,\  r^*}^{k} r_j(\lambda_j -  x_{n-1}) \ge  r^*D - D^{(m)}
 \sum_{s=1}^{l(m)} r_{j_s} -   D^{(m+1)} \sum_{s=1}^{l(m+1)} r_{l_s} $$
$$- \left(n-r_1-  r^*- r_*-  \sum_{s=1}^{l(m)} r_{j_s}-  \sum_{s=1}^{l(m+1)} r_{l_s}\right)D.$$
But $x_{n-1}- \lambda_*= \hbox{span}(f)- D$.  Therefore,
$$r_*  \hbox{span}(f) + \left(n-r_1- 2(  r^*+ r_*)\right) D \ge (D- D^{(m)})  \sum_{s=1}^{l(m)} r_{j_s} +
(D- D^{(m+1)})  \sum_{s=1}^{l(m+1)} r_{l_s}.$$
The right-hand side of the latter inequality is, obviously,  greater or equal to $r_0 \left(l(m)+ l(m+1)\right) (D- D^{(m)}) $, where $1 \le r_0= \hbox{min}_{1\le j\le k} (r_j)$.  Moreover,  since $ \hbox{span}(f) \le 2D$,  the left-hand side does not exceed $\left(n-r_1\right) D-  r^* \hbox{span}(f)$. Thus  we come up with the inequality
$$r_0 \left(l(m)+ l(m+1)\right) (D- D^{(m)}) \le \left(n-r_1\right) D-  r^* \hbox{span}(f)$$
or since $D- D^{(m)} > 0$ \ ($m \ge r$), it becomes
$$l(m)+ l(m+1)\le \frac{\left(n-r_1\right) D-  r^* \hbox{span}(f)}{r_0(D- D^{(m)}) }.\eqno(43)$$
Meanwhile, appealing to (16), we get similarly 
$$n(n-1) ( x_{n-1}- x_{n-2})^2 = r^*D^2 + r_* (\lambda_*- x_{n-1})^2 +  \sum_{j=2, \ r_j\neq r_*,\  r^*}^{k}
 r_{j}(\lambda_j-  x_{n-1})^2 $$$$
\le  r^*D^2 + r_*\left(\hbox{span}(f)- D\right)^2+  \left[D^{(m)}\right]^2 \sum_{s=1}^{l(m)} r_{j_s} +   \left[D^{(m+1)}\right]^2 \sum_{s=1}^{l(m+1)} r_{l_s}$$$$+ \left(n-r_1-  r^*- r_*-  \sum_{s=1}^{l(m)} r_{j_s}-
 \sum_{s=1}^{l(m+1)} r_{l_s}\right)D^2.$$
 Therefore, analogously to (43), we arrive at the inequality
 $$l(m)+ l(m+1)\le   \frac{(n-r_1)D^2 + r_* \left[\hbox{span}(f)\right]^2 - n(n-1) ( x_{n-1}- x_{n-2})^2 -2D r_*\  \hbox{span}(f)}{r_0(D^2- \left[D^{(m)}\right]^2) }.$$
{\bf Proposition 3}. {\it There exists no CA- polynomial with only real roots of degree $n$ such that }
$$  \hbox{span}(f) >  \left(r^*\right)^{-1} \left[ (n-r_1-r_0)D + r_0 D^{(m)}\right],\   m\ge r.\eqno(44)$$
\begin{proof}  Under condition (44), the right-hand side of (43) is less than one.  Thus $l(m)= l(m+1)=0$ and Lemma 11 completes the proof.
\end{proof}
Let $m=n-2$.   Then since $l(n-1)=0$, inequality (43) becomes
$$l(n-2) \le  \frac{ (n-r_1)D -  r^* \hbox{span}(f)}{r_0(D-   \left|x_{n-1}- x_{n-2}\right|)}.\eqno(45)$$
{\bf Proposition 4}. {\it There exists no CA- polynomial with only real roots of degree $n$ such that }
$$D <  \left[r^* \  \sqrt{ \frac{n^2-r_1}{n-r_1-r_2}} -r_0\right] \frac{ \left|x_{n-1}- x_{n-2}\right|}{n-r_1-r_0}.\eqno(46)$$

\begin{proof}  Indeed, employing the lower bound (38) for $ \hbox{span}(f)$, we find that under condition (46) the right-hand side of (45) is strictly less than one.  Consequently,  $l(n-2)=0$ and owing to Corollary 1 $f$ is trivial. If the maximum of multiplicities $r > n-2$,    $f$ has at most 2 distinct zeros and  it is trivial via Corollary 3.
\end{proof}

Finally, we prove 

{\bf Proposition 5}. {\it Let  CA- polynomial with only real roots exist.  Then it has the property
$$  \frac{d}{D} \le \sqrt{\frac{2(n-m-1)}{2(k-1)-1}},\eqno(47)$$
where $d, D$ are defined by $(34), (35)$, respectively,  and $m,\ m+1$ belong to the interval $\left[r, \   {1\over 2}\left(1-{1\over r_0}\right)(n-1) \right)$.}
\begin{proof}  Since $m,\ m+1$ are chosen from  the interval $\left[r, \   {1\over 2}\left(1-{1\over r_0}\right)(n-1) \right)$,  condition (24) holds for these values. Hence assuming the existence of the CA-polynomial, we return to the Sz.-Nagy type identity (25) to  have the estimate 
$$0 \ge  l(m)  \left( r_{0} {(n-m)(n-m -1)\over n(n-1)} -1\right) d^2 +  \left(k-1-l(m)\right) d^2- (n-m-l(m)) D^2$$
$$\ge  (k-1) d^2- (n-m) D^2  + l(m) (D^2-d^2).$$ 
Writing the same inequality for $m+1$
$$0 \ge  (k-1) d^2- (n-m-1) D^2  + l(m+1) (D^2-d^2)$$ 
and adding two inequalities, we find
$$0 \ge  2 (k-1) d^2- (2(n-m) -1) D^2  +(l(m)+  l(m+1)) (D^2-d^2),$$
which means
$$l(m)+  l(m+1) \le \frac{(2(n-m) -1) D^2 - 2 (k-1) d^2}{D^2-d^2}.$$
So, for the existence of the CA-polynomial it is necessary that the right-hand side of the latter inequality is greater than  or equal to 1.  Thus we come up  with condition (47) and complete the proof. 
\end{proof}

{\bf  Acknowledgment}.   The present investigation was supported, in
part, by the "Centro de Matem{\'a}tica" of the University of Porto.

\bibliographystyle{amsplain}

\end{document}